\documentclass{amsart}
\usepackage{amssymb}
\usepackage{amsmath}

\newtheorem{theorem}{Theorem}[section]

\newtheorem{proposition}[theorem]{Proposition}

\theoremstyle{definition}\newtheorem{remark}{Remark}\newtheorem{question}{Question}

\newcommand{\euO}{\mathfrak O}
\newcommand{\euP}{\mathfrak P}
\newcommand{\euD}{\mathfrak D}

\begin{document}
\title[Valuation criterion]{A valuation criterion for normal basis generators in local fields of characteristic $p$}

\author{G. Griffith Elder}
\email{elder@unomaha.edu} \address{Department of Mathematics \\ University of Nebraska at Omaha\\ Omaha, NE 68182-0243 U.S.A.}

\date{February 12, 2008}
%\classno{11S15}

%\maketitle

\begin{abstract}
Let $K$ be a complete local field of characteristic $p$ with perfect
residue field.  Let $L/K$ be a finite, fully ramified, Galois
$p$-extension.  If $\pi_L\in L$ is a prime element, and $p'(x)$ is the
derivative of $\pi_L$'s minimal polynomial over $K$, then the relative
different $\euD_{L/K}$ is generated by $p'(\pi_L)\in L$.  Let $v_L$ be
the normalized valuation normalized with $v_L(L)=\mathbb{Z}$.
We show that any element $\rho\in L$ with
$v_L(\rho)\equiv -v_L(p'(\pi_L))-1\bmod[L:K]$ generates a normal
basis, $K[\mbox{Gal}(L/K)]\cdot\rho=L$.  This criterion is tight:
Given any integer $i$ such that $i\not\equiv
-v_L(p'(\pi_L))-1\bmod[L:K]$, there is a $\rho_i\in L$ with
$v_L(\rho_i)=i$ such that $K[\mbox{Gal}(L/K)]\cdot\rho_i\subsetneq L$.
\end{abstract}

\maketitle

The Normal Basis Theorem states that in a finite Galois extension
$L/K$ with $G=\mbox{Gal}(L/K)$, there is an element $\rho\in L$, called a
normal basis generator, whose conjugates $\{\sigma\rho:\sigma\in G\}$
provide a basis for $L$ over $K$.  In the setting of local field
extensions, the most important property of an element is its
valuation, and so \cite{elder:blms} asked whether there is a valuation
criterion: Is there a valuation (an integer certificate) that
guarantees that any element bearing this valuation is a normal basis
generator?  

Let $K$ be a complete local field with a perfect residue field of
characteristic $p$. So the characteristic of $K$ is $0$ or $p$. 
It is
not too hard to see that as a necessary condition, if a valuation
criterion exists for a finite Galois extension $L/K$, then $L/K$ must be
fully ramified and have order a power of $p$. But are these necessary
conditions also sufficient?  In this paper we give an affirmative
answer to this question in the case where $K$ has characteristic $p$.

For $K$ of characteristic $0$ and regular (so $K$ does not contain the
$p$th roots of unity), a valuation criterion for fully ramified
elementary abelian $p$-extensions is given in \cite{elder:blms}. For
$K$ of characteristic $p$, a valuation criterion for fully ramified
abelian $p$-extensions is given in \cite{thomas}.  In both cases, the
valuation criterion is described in terms of the largest ramification
break number associated with $L/K$. The main contribution of this
paper is a restatement of that criterion in terms of the exponent of
the relative different (see remark following proposition). The
relative different satisfies $\euD_{L/K}=(p'(\pi_L))$ with $\pi_L$ a
prime element of $L$ and $p(x)$ the minimal polynomial of $\pi_L$ over
$K$ \cite[III \S6 Lemma 2 Corollary 2]{serre:local}. This means that
we can state our main result as follows:

%\begin{theorem*}
\begin{theorem}
Let $K$ be a complete local field of characteristic $p$ with perfect residue
field. Let $L/K$ be a finite, Galois extension with $G=\mbox{\rm
  Gal}(L/K)$, and let $v_L$ be the valuation normalized so that
$v_L(L)=\mathbb{Z}$.  If $L/K$ is fully ramified and $[L:K]=p^n$ for
some integer $n$, there is a valuation criterion: Let $\pi_L\in L$ be
a prime element and let $p(x)$ be its
minimal polynomial over $K$. Then any element $\rho\in L$ with
$v_L(\rho)\equiv -v_L(p'(\pi_L))-1\bmod[L:K]$ generates a normal basis
for $L/K$. So $L=K[G]\cdot\rho$.

Moreover the assumption that $L/K$ is a fully ramified $p$-extension
is necessary and the criterion under that assumption is
tight: Outside of the assumption, given any $i\in\mathbb{Z}$, or under
the assumption, given any $i\not\equiv -v_L(p'(\pi_L))-1\bmod[L:K]$,
there is a $\rho_i\in L$ with $v_L(\rho_i)=i$ such that $L\supsetneq
K[G]\cdot\rho_i$.
%\end{theorem*}
\end{theorem}

\begin{proof} Recall the definition of the 
ramification groups $G_i$ \cite[IV]{serre:local}. We begin by assuming
that $G=G_1$, which is equivalent to the assumption that $L/K$ is a
fully ramified extension of degree $p^n$ for some $n$.  Let
$d=v_L(p'(\pi_L))$ and let $\rho\in L$ with $v_L(\rho)\equiv -d-1\bmod
p^n$.  Express $\rho$ in terms of the field basis
$B=\{\pi_L^i/p'(\pi_L):i=0,\ldots ,p^n-1\}$ for $L/K$. It is a result
of Euler that
$$\mbox{Tr}_G\frac{\pi_L^i}{p'(\pi_L)}=\begin{cases}0&\mbox{for }0\leq
i\leq p^n-2\\ 1&\mbox{for }i=p^n-1\end{cases}$$ \cite[III \S6 Lemma
  2]{serre:local}.  Since $v_L(\rho)\equiv
v_L(\pi_L^{p^n-1}/p'(\pi_L))\bmod p^n$, the coefficient of
$\pi_L^{p^n-1}/p'(\pi_L)$ in the expression for $\rho$ must be
nonzero. Therefore $\mbox{Tr}_G\rho\neq 0$. Moreover, we can replace
$B$ by the alternate basis $\{\rho,\pi_L^i/p'(\pi_L):i=0,\ldots
,p^n-2\}$, and get
$$L=K\cdot\rho+\sum_{i=0}^{p^n-2}K\cdot \frac{\pi_L^i}{p'(\pi_L)}.$$

The Normal Basis Theorem, stated in terms of Tate cohomology, says
that $\hat{H}^{-1}(G,L)=0$ \cite[VIII \S1 \& X \S1 Proposition
  1]{serre:local}.  This means that any element $\eta\in L$ with
$\mbox{Tr}_G\eta=0$ satisfies $\eta\in \mathcal{I}_G\cdot L$, where
$\mathcal{I}_G=(\sigma-1:\sigma\in G)$ is the augmentation ideal of
the group ring $K[G]$.  Therefore $L=K\cdot\rho +\mathcal{I}_G\cdot
L$, but also
$$L=K[G]\cdot \rho +\mathcal{I}_G\cdot L.$$ Now notice that because
$G$ is a $p$-group and $K$ has characteristic $p$, $\mathcal{I}_G$ is
also the Jacobson radical of $K[G]$.  Thus, by Nakayama's Lemma,
$L=K[G]\cdot\rho$. We have proven the criterion.  The elements
$\pi_K^m\pi_L^{p^n-1}/p'(\pi_L)$ for $i\neq p^n-1$ show that the
criterion is sharp.  (All of this argument applies equally well in
characteristic $0$, except for one step: In characteristic $0$, the
augmentation ideal $\mathcal{I}_G$ is not the Jacobson radical of
$K[G]$.)

Now we address $G\supsetneq G_1$. Let $L'=L^{G_1}$ and
$L''=L^{G_0}$. So $L''/K$ is unramified and $L'/L''$ is fully but
tamely ramified.  Replace $K$ and $G$, in our earlier argument with
$L'$ and $G_1$ respectively.  So $d$ is defined in terms of the
relative different of the fully ramified $p$-part of the extension,
$\euP_L^d=\euD_{L/L'}$ and $p^n=[L:L']$.  From earlier work,
given any integer $i\not\equiv -d-1\bmod p^n$ there is a $\rho\in L$ with
$v_L(\rho)=i$ such that $\mbox{Tr}_{G_1}\rho=0$. So
$K[G]\rho\subsetneq L$.

To consider the case $i\equiv -d-1 \bmod p^n$, notice that the trace
$\mbox{Tr}_{G_1}$ maps fractional ideals of $\euO_L$, the ring of
integers of $L$, to fractional ideals of $\euO_{L''}$. Indeed, using
basis $B$ from above,
$\mbox{Tr}_{G_1}\euP_L^{kp^n-d-1}=\euP_{L'}^{k-1}$. Moreover
$\mbox{Tr}_{G_1}\euP_L^{kp^n-d}\subseteq \euP_{L'}^{k}$, because of
\cite[III \S3 Proposition 7]{serre:local}. Observe that this means
that given any element $\tau\in L'$ with $v_{L'}(\tau)=k-1$, there is
a $\rho_{\tau}\in L$ with $v_L(\rho_{\tau})=kp^n-d-1$ such that
$\mbox{Tr}_{G_1}(\rho_{\tau})=\tau$.

To use this observation notice that because $L'/L''$ is tamely
ramified, there is a prime element of $L''$, namely $\pi_{L''}$, such
that $L'=L''(\sqrt[e]{\pi_{L''}})$ with $p\nmid e$ \cite[II \S3.5
  Proposition]{fesenko}.  So for $k\not\equiv 1\bmod e$, let $\tau=
\sqrt[e]{\pi_{L''}}^{k-1}$. Since
$\mbox{Tr}_{G_0/G_1}\sqrt[e]{\pi_{L''}}^{k-1}=0$, $\mbox{Tr}_{G_0}\rho
=0$ and $K[G]\rho\subsetneq L$. For $k\equiv 1\bmod e$, let
$\tau=\pi_K^{(k-1)/e}$ and let $\bar{\sigma}$ be any non-trivial
element in $G/G_1$.  Then $(\sigma-1)\pi_K^{k-1}=0$. Thus
$(\sigma-1)\mbox{Tr}_{G_1}\rho=0$ and $K[G]\rho\subsetneq L$.
\end{proof}

Now we connect the valuation criterion of this paper with that of 
\cite{elder:blms}, \cite{thomas}.

\begin{proposition}
%\begin{proposition*}
Let $K$ be a complete local field (characteristic $0$ or $p$) with perfect
residue field of characteristic $p$.  Let $L/K$ be a finite, fully
ramified, Galois $p$-extension. Let $\pi_L\in L$ be
a prime element and let $p(x)$ be its
minimal polynomial over $K$. Then
$$v_L(p'(\pi_L))+1\equiv
p^nu_m-b_m\bmod[L:K]$$ where $b_m, u_m$ are the
largest ramification break numbers in lower and upper numbering
respectively. 
%\end{proposition*}
\end{proposition}
%\begin{remark*}
\begin{remark}
When $L/K$ is abelian $u_m\in\mathbb{Z}$ by the Hasse-Arf Theorem and
so we get the valuation criterion of \cite{elder:blms}, \cite{thomas}, namely
$v_L(\rho)\equiv b_m\bmod p^n$.
%\end{remark*}
\end{remark}
\begin{proof}
Let $d=v_L(p'(\pi_L))$ and recall that $\euD_{L/K}=(p'(\pi_L))$.
Let $G=G_1$ and $[L:K]=p^n$.
Recall the ramification filtration $G=G_1\supseteq G_2\supseteq\cdots
$ where $G_i=\{\sigma\in G:v_L((\sigma-1)\pi_L)\geq i+1\}$ \cite[IV
  \S1]{serre:local}. The break numbers (in lower numbering) are those
integers $i$ such that $G_i\supsetneq G_{i+1}$. Let $b_1<b_2<\cdots
<b_m$ be the list of break numbers with $b_m$ being the maximal break
(again in lower numbering). Let $g_i=|G_{b_i}|$, the number of
elements in $G_{b_i}$. Then by \cite[IV \S1 Proposition
  4]{serre:local},
$$d=(1+b_1)(g_1-1)+\sum_{i=2}^m (b_i-b_{i-1})(g_i-1)
=(1+b_1)g_1-b_m-1+\sum_{i=2}^m (b_i-b_{i-1})g_i .$$ Moreover, we can
convert the lower numbering to upper numbering using the Herbrand
function $\varphi$ \cite[IV \S3]{serre:local}. Be careful to notice a
small difference in notation: We use $g_i=|G_{b_i}|$ whereas \cite[IV
\S3]{serre:local} uses $g_{b_i}=|G_{b_i}|$. The largest break
number in upper numbering is therefore
$$u_m=\varphi(b_m)=\frac{1}{p^n}\left (b_1g_1+\sum_{i=2}^m
(b_i-b_{i-1})g_i\right ).$$ Thus $d+1=g_1-b_m+p^nu_m$, where
$g_1=p^n$.
\end{proof}

We end the paper with a natural 
%\begin{question*}
\begin{question}
Does the statement of the Theorem (modified appropriately to address
\cite[Example 1]{elder:blms}) also hold in characteristic zero?
%\end{question*}
\end{question}

%\begin{acknowledgements}\label{ackref} The author thanks Nigel Byott and Bart de Smit for their suggestions.\end{acknowledgements}

%\affiliationone{G.~Griffith Elder\\ Department of Mathematics\\ University of Nebraska at Omaha\\ Omaha, NE 68182-0243 U.S.A. \email{elder@unomaha.edu}}


\begin{thebibliography}{9}

\bibitem{desmit}
%{\bibname B. de Smit \and L. Thomas},
{B. de Smit \and L. Thomas},
`Local {G}alois module structure in positive characteristic and continued fractions',
{\em Arch. Math. (Basel)} 88 (2007), no. 3, 207--219.

\bibitem{elder:blms}
%{\bibname N. P. Byott \and G. G. Elder},
{N. P. Byott \and G. G. Elder},
`A valuation criterion for normal bases in elementary abelian extensions',
{\em Bull. Lond. Math. Soc.} 39 (2007), no. 5, 705--708.

\bibitem{fesenko}
%{\bibname I. B. Fesenko \and S. V. Vostokov},
{I. B. Fesenko \and S. V. Vostokov},
`Local fields and their extensions', (American Mathematical Society, Providence RI, 2002).

\bibitem{serre:local} 
%{\bibname J.-P. Serre}, 
{J.-P. Serre}, 
`Local fields',
  (Springer-Verlag, New York, 1979).

\bibitem{thomas} 
%{\bibname L. Thomas}, 
{L. Thomas}, 
`A valuation criterion for
  normal basis generators in equal positive characteristic', preprint:
  August 9, 2007.

\end{thebibliography}
\end{document}